\documentclass[a4paper,reqno,twoside,draft]{amsart}

\usepackage{latexsym,amssymb}

\theoremstyle{plain}
\newtheorem{propn}{Proposition}[section]
\newtheorem{thm}[propn]{Theorem}
\newtheorem{lemma}[propn]{Lemma}
\newtheorem{cor}[propn]{Corollary}
\newtheorem*{thm*}{Theorem}

\theoremstyle{definition}

\theoremstyle{remark}
\newtheorem*{rem}{Remark}
\newtheorem*{rems}{Remarks}

\newcommand{\Hil}{\mathsf{H}}
\newcommand{\hil}{\mathsf{h}}
\newcommand{\Kil}{\mathsf{K}}

\newcommand{\Fock}{\mathcal{F}}
\newcommand{\Exp}{\mathcal{E}}
\newcommand{\e}[1]{\varepsilon (#1)}
\newcommand{\w}[1]{\varpi (#1)}
\newcommand{\step}{\mathbb{S}}

\newcommand{\init}{\mathfrak{h}}
\newcommand{\D}{\mathfrak{D}}
\newcommand{\Dvd}{\D^{V,d}}
\newcommand{\Dd}{\D^d}
\newcommand{\DvT}{\D^{V,\iT}}
\newcommand{\DvD}{\D^{V,\iD}}
\newcommand{\Dvo}{\D^{V,0}}

\newcommand{\noise}{\mathsf{k}}
\newcommand{\khat}{\wh{\noise}}
\newcommand{\iT}{\mathsf{T}}
\newcommand{\iD}{\mathsf{D}}
\newcommand{\iDhat}{{\wh{\iD}}}

\newcommand{\Vtilde}{{\widetilde{V}}}

\newcommand{\Q}{\mathcal{Q}}

\newcommand{\Qcd}{Q^{c,d}}
\newcommand{\Qab}{Q^{a,b}}
\newcommand{\VQcd}{{}^V\!\!\Qcd}
\newcommand{\Qncd}{{}^{(n)}\!\Qcd}

\newcommand{\Gcd}{G_{c,d}}
\newcommand{\GVtildecd}{G^\Vtilde_{c,d}}
\newcommand{\Gh}[2]{G_{({#1},{#2})}}

\newcommand{\VPcd}{{}^V\!\!P^{c,d}}
\newcommand{\VPad}{{}^V\!\!P^{a,d}}
\newcommand{\VPcb}{{}^V\!\!P^{c,b}}
\newcommand{\VPab}{{}^V\!\!P^{a,b}}

\newcommand{\HVcd}{H^V_{c,d}}
\newcommand{\HVtildecd}{H^\Vtilde_{c,d}}

\newcommand{\tot}{_{[0,t[}}
\newcommand{\tT}{_{[t,T[}}
\newcommand{\toT}{_{[0,T[}}
\newcommand{\fromt}{_{[t,\infty[}}

\newcommand{\Form}{\mathbb{F}}
\newcommand{\Op}{\mathcal{O}}
\newcommand{\Mat}{\mathrm{M}}
\newcommand{\bfc}{\mathbf{c}}

\newcommand{\Fab}{F^\al_\be}

\newcommand{\al}{\alpha}
\newcommand{\be}{\beta}

\newcommand{\Ga}{\Gamma}
\newcommand{\si}{\sigma}

\newcommand{\bvarpi}{\boldsymbol{\varpi}}

\newcommand{\Real}{\mathbb{R}}
\newcommand{\Rplus}{\Real_+}
\newcommand{\Comp}{\mathbb{C}}
\newcommand{\Nat}{\mathbb{N}}
\newcommand{\Int}{\mathbb{Z}}

\newcommand{\ip}[2]{\langle #1, #2 \rangle}
\newcommand{\bip}[3][\big]{#1\langle #2, #3 #1\rangle}
\newcommand{\norm}[1]{\lVert #1 \rVert}
\newenvironment{sbmatrix}
{\bigl[\begin{smallmatrix}}{\end{smallmatrix}\bigr]}

\newcommand{\wh}{\widehat}
\newcommand{\wt}{\widetilde}
\newcommand{\ol}{\overline}

\newcommand{\schur}{\bullet}
\newcommand{\ot}{\otimes}

\newcommand{\op}{\oplus}
\newcommand{\uot}{\underline{\otimes}}

\newcommand{\To}{\rightarrow}
\newcommand{\Tends}{\rightarrow}
\newcommand{\Implies}{\Rightarrow}

\newcommand{\loc}{\text{\tu{loc}}}

\newcommand{\ti}{\textit}
\newcommand{\tu}{\textup}
\newcommand{\dfn}{\ti}

\DeclareMathOperator{\Dom}{Dom}
\DeclareMathOperator{\Ran}{Ran}
\DeclareMathOperator{\Lin}{Lin}
\DeclareMathOperator{\Aff}{Aff}

\DeclareMathOperator{\re}{Re}

\DeclareMathOperator{\diag}{diag}
\DeclareMathOperator{\wlim}{\mbox{w-lim}}

\newenvironment{alist}
{

\begin{enumerate}}
{\end{enumerate}}

\newenvironment{rlist}
{

\begin{enumerate}}
{\end{enumerate}}

\numberwithin{equation}{section}

\begin{document}

\title[Quantum Stochastic Operator Cocycles]{Quantum Stochastic
Operator Cocycles
\\ via Associated Semigroups}
\author{J.~Martin Lindsay}
\address{Department of Mathematics and Statistics \\ Lancaster
University \\ Lancaster LA1 4YF \\ UK }
\email{j.m.lindsay@lancaster.ac.uk}
\author{Stephen J.~Wills}
\address{School of Mathematical Sciences \\ University College Cork \\
Cork \\ Ireland}
\email{s.wills@ucc.ie}

\begin{abstract}
A recent characterisation of Fock-adapted contraction operator
stochastic cocycles on a Hilbert space, in terms of their associated
semigroups, yields a general principle for the construction of such
cocycles by approximation of their stochastic generators. This leads
to new existence results for quantum stochastic differential
equations. We also give necessary and sufficient conditions for a
cocycle to satisfy such an equation.
\end{abstract}

\maketitle

\section{Introduction}

In this paper we study the functional equation
\[
V_0 = I, \quad V_{r+t} = V_r \si_r (V_t) \ \text{ for all } r,t
\geqslant 0
\]
for a family of contractions on $\init \ot \Fock$ adapted to the Fock
operator filtration. Here $\Fock$ is the symmetric Fock space over
$L^2(\Rplus; \noise)$, $\init$ and $\noise$ are fixed but arbitrary
Hilbert spaces, and $(\si_t)_{t \geqslant 0}$ is the endomorphism
semigroup of shifts, ampliated to $B(\init \ot \Fock)$. We call such a
family a \dfn{left contraction cocycle on $\init$} with noise
dimension space $\noise$.

Contraction cocycles may be constructed by solving quantum stochastic
differential equations of Hudson-Parthasarathy type. By means of a
recent characterisation of such cocycles, in terms of an associated
family of semigroups (Theorem~\ref{recon thm}), we provide a new
method of constructing cocycles which in turn leads to new existence
results for QSDEs. When the driving noise is infinite dimensional the
coefficient of a QSDE is naturally given as a sesquilinear
operator-valued map or, in terms of a coordinate system for the noise
dimension space $\noise$, as an infinite matrix $[\Fab]$. We show that
if a process satisfying such a \emph{form QSDE} is contractive and
strongly measurable then the coefficient is necessarily given by an
\emph{operator}, equivalently the matrix must be \emph{semiregular}.
We also give necessary and sufficient conditions, of weak
differentiability type, for a strongly continuous left contraction
cocycle to satisfy a QSDE.

This paper builds on work of Accardi, Fagnola, Journ\'{e}, Mohari 
and the authors (\cite{FF8, journe, Moh, AJL, RHP}), extending known 
results for Markov-regular cocycles and QSDEs with bounded coefficients
(\cite{HuP, HuL, mother, father}; see \cite{meyer, L} and references
therein). Our development of the theory is coordinate-free, moreover a
technical feature of the work is that no separability assumptions are
imposed on either the initial space $\init$ or the noise dimension
space $\noise$. This freedom is useful for certain applications such
as the stochastic dilation of quantum dynamical semigroups
(\cite{modules}). A different approach to the characterisation and
construction of cocycles through semigroup methods has been outlined
by Liebscher (\cite{Liebscher}).

\subsection{General notations}

The algebraic tensor product is denoted $\uot$, with $\ot$ reserved
for the tensor product of Hilbert spaces and their operators. For a
vector $\xi$ in a Hilbert space $\Kil$, operators $E_\xi: \Hil \To
\Hil \ot \Kil$ and $E^\xi: \Hil \ot \Kil \To \Hil$ are defined by
\[
E_\xi u = u \ot \xi \quad \text{ and } \quad E^\xi =
(E_\xi)^*;
\]
with context indicating the Hilbert space $\Hil$, and moreover the
elementary tensor $u \ot \xi$ is usually abbreviated to $u\xi$. Note
that $\xi \mapsto E_\xi$ is an isometry. For Hilbert spaces $\Hil $
and $\Hil'$ and a dense subspace $\mathcal{D}$ of $\Hil$, $\Op
(\mathcal{D}; \Hil')$ denotes the linear space of operators $\Hil \To
\Hil'$ with domain $\mathcal{D}$; $\Op (\mathcal{D})$ abbreviates $\Op
(\mathcal{D}; \Hil)$. For $f \in L^2 (\Rplus; \noise)$ and $I \subset
\Rplus$, $f_I$ denotes the function that agrees with $f$ on $I$ and is
zero elsewhere, and $c_I$ denotes the function equal to $c$ on $I$ and
zero elsewhere, when $c$ is a vector in $\noise$.

\subsection{Fock space}

We use normalised exponential vectors in $\Fock = \Ga \bigl( L^2
(\Rplus; \noise) \bigr)$, the symmetric Fock space over the Hilbert
space $L^2 (\Rplus; \noise)$. These are defined by $\w{f}:=
\norm{\e{f}}^{-1} \e{f}$ where $ \e{f} = (1,f,(2!)^{-1/2} f^{\ot 2},
\ldots)$ for $f \in L^2(\Rplus; \noise)$. The function
\begin{equation} \label{chi}
\chi: \noise \times \noise \To \Comp, \quad (c,d) \mapsto \tfrac{1}{2}
\bigl( \norm{c}^2 + \norm{d}^2 \bigr) - \ip{c}{d},
\end{equation}
which governs their inner product: $\ip{\w{f}}{\w{g}} = \exp
\bigl(-\int \chi (f(s), g(s)) \, ds \bigr)$, also plays a role. The
subspace $\Exp(S): = \Lin \{\e{f}: f \in S \}$ is dense in $\Fock$ for
various useful subsets of $L^2 (\Rplus; \noise)$, for example
\[
\step_\iT:= \{ f \in \step: f \text{ is } \iT\text{-valued}\}
\]
where $\step := \Lin \{c\tot: c \in \noise, t>0\}$, and $\iT$ is a
total subset of $\noise$ containing $0$; we write $\Exp_\iT$ for $\Exp
(\step_\iT)$. Examples of such sets $\iT$ include (not necessarily
normalised) orthogonal bases, augmented by $0$.

We shall need a refinement of the basic estimate
\begin{equation} \label{basic}
\norm{\w{f} -\w{g}} \leqslant \norm{\e{f} -\e{g}} \leqslant \norm{f-g}
e^{\frac{1}{2} (\norm{f} + \norm{g})^2}
\end{equation}
obtained by viewing $L^2(\Real_+; \noise)$ as a subspace of $\Fock$,
namely
\begin{equation} \label{refined}
\norm{\e{f} -\e{g} -(f-g)} \leqslant \norm{f-g} \bigl( \norm{f}
+\norm{g} \bigr) e^{\frac{1}{2} (\norm{f} + \norm{g})^2}.
\end{equation}

Letting $\Fock_I$ denote the symmetric Fock space over $L^2 (I;
\noise)$, for a subinterval $I$ of $\Rplus$, the tensor factorisation 
$\Fock \cong \Fock\tot \ot \Fock\fromt$ is given by continuous linear 
extension of the correspondence
\[
\w{f} \longleftrightarrow \w{f|\tot} \w{f|\fromt}.
\]

\subsection{Operator processes}

A family of operators $(X_t)_{t \geqslant 0}$ in $ B(\init \ot \Fock)$
is \dfn{adapted} if
\[
X_t \in B(\init \ot \Fock\tot) \ot I_{\Fock\fromt} \ \text{ for all }
t,
\]
that is, if each $X_t$ only acts nontrivially before time $t$, and is
called an (\dfn{operator}) \dfn{process} if furthermore it is weakly
measurable (i.e.\ if $t \mapsto \ip{\xi}{X_t \zeta}$ is measurable for
all $\xi, \zeta \in \init \ot \Fock$). Here we are concerned with
contraction operator-valued processes, which we refer to simply as
\dfn{contraction processes} on $\init$.

The right shift $s_t$ and time reversal map $r_t$ on $L^2(\Rplus;
\noise)$ are
\[
(s_t f)(u) = \begin{cases} 0 & \text{if } u < t, \\ f(u-t) & \text{if
} u \geqslant t, \end{cases}
\quad \text{ and } \quad
(r_t f)(u) = \begin{cases} f(t-u) & \text{if } u \leqslant t, \\ f(u)
& \text{if } u > t. \end{cases}
\]
Their second quantisations ampliated to $\init \ot \Fock$ are the
isometry $S_t$ and self-adjoint unitary $R_t$ respectively, given by
\begin{equation} \label{S and R}
S_t u\e{f} = u\e{s_t f}, \quad R_t u\e{f} = u\e{r_t f}.
\end{equation}
Thinking of $s_t$ as a unitary $L^2(\Rplus; \noise) \To
L^2([t,\infty[; \noise)$ we get the Hilbert space isomorphism $\Fock
\cong \Fock\fromt$, and the algebra isomorphism $B(\init \ot \Fock)
\cong B(\init \ot \Fock\fromt)$. The latter algebra is viewed as a
subalgebra of the former:
\[
B(\init \ot \Fock) \cong B(\init) \ot I_{\Fock\tot} \ot
B(\Fock\fromt) \subset B(\init \ot \Fock).
\]
Then, for $Y \in B(\init \ot \Fock)$, $\si_t(Y) \in B(\init \ot
\Fock)$ denotes the result of carrying out these identifications; more
concretely it is determined by the identity
\[
\ip{u\e{f}}{\si_t (Y) v\e{g}} = \ip{u\e{s^*_t f}}{Y v\e{s^*_t g}}
\ip{\e{f\tot}}{\e{g\tot}}.
\]
The family $(\si_t)_{t \geqslant 0}$ is a pointwise weakly continuous
semigroup of normal endomorphisms of $B(\init \ot \Fock)$.

\subsection{Perturbation}

We end this introduction by quoting a dissipative generalisation of
the Kato-Rellich Theorem whose symmetric form is well-suited to our
purposes. Recall that a $C_0$-semigroup is contractive if and only if
its generator is dissipative.

\begin{thm} \label{Gustafson}
Let $A$ and $B$ be densely defined dissipative operators on a Banach
space with the same domain $\mathcal{D}$, and suppose that there are
constants $\lambda, \mu \geqslant 0$ with $\lambda <1$ such that their
difference $D=A-B$ satisfies
\begin{equation} \label{A-B}
\norm{Dv} \leqslant \lambda \bigl( \norm{Av} + \norm{Bv} \bigr) + \mu
\norm{v}, \quad v \in \mathcal{D}.
\end{equation}
Then $\Dom \ol{A} = \Dom \ol{B}$, moreover $\ol{A}$ is a
$C_0$-semigroup generator if and only if $\ol{B}$ is.
\end{thm}

This result is due to Gustafson (see~\cite{ReS}, Theorem~X.50). Note
that if~\eqref{A-B} holds then $D$ is $A$-bounded with relative bound
at most $2\lambda /(1-\lambda)$.

\section{Cocycles and Semigroups}

Left contraction cocycles on $\init$ have been defined in the
introduction. An adapted family of contractions $U=(U_t)_{t \geqslant
0}$ on $\init \ot \Fock$, satisfying $U_0 = I$ and $U_{r+t} = \sigma_r
(U_t) U_r$ for $r,t \geqslant 0$, is called a \dfn{right contraction
cocycle}. Thus $U$ is a right contraction cocycle if and only if $U^*
:= (U^*_t)_{t \geqslant 0}$ is a left contraction cocycle.

\subsection{The semigroup decomposition}

For a contraction process $V = (V_t)_{t \geqslant 0}$ on $\init$,
define the following operators on $\init$:
\[
\VQcd_t = E^{\w{c\tot}} V_t E_{\w{d\tot}}, \quad c,d \in \noise, t
\geqslant 0.
\]
These `sliced' operators allow one to determine whether or not a
process $V$ is a left contraction cocycle. Note that they are all
contractions.

\begin{propn}[\cite{father}] \label{Semi decomp}
Let $V$ be a contraction process on $\init$, and let $\iT$ and
$\iT^\dagger$ be any total subsets of $\noise$ containing $0$. Then
the following are equivalent\tu{:}
\begin{rlist}
\item
$V$ is a left contraction cocycle.
\item
For each choice of $c \in \iT^\dagger$ and $d \in \iT$, $(\Qcd_t :=
\VQcd_t)_{t\geqslant 0}$ is a contraction semigroup on $\init$, and
for all $f \in \step_{\iT^\dagger}$ and $g \in \step_\iT$
\begin{equation} \label{semi decomp}
E^{\w{f\tot}} V_t E_{\w{g\tot}} = Q^{f(t_0),g(t_0)}_{t_1-t_0} \cdots
Q^{f(t_n),g(t_n)}_{t-t_n}
\end{equation}
where $\{0=t_0 \leqslant t_1 \leqslant \cdots \leqslant t_n \leqslant
t\}$ contains the discontinuities of $f\tot$ and $g\tot$, and
\emph{right-continuous versions} are used for the evaluations.
\end{rlist}
\end{propn}
We refer to $\{\VQcd: c,d \in \noise\}$ as the cocycle's
\dfn{associated semigroups}.

\begin{rem}
The same holds for right contraction cocycles except that the product
in~\eqref{semi decomp} is in the reverse order. It follows that
$(U_t)_{t \geqslant 0}$ is a right contraction cocycle if and only if
$(R_t U_t R_t)_{t \geqslant 0}$ defines a left contraction cocycle,
where the operators $R_t$ are defined in~\eqref{S and R}.
\end{rem}

For a left contraction cocycle $V$, we refer to the left contraction
cocycle defined by $(R_t V_t^* R_t)_{t \geqslant 0}$ as the
(\dfn{Journ\'{e}}) \dfn{dual} of $V$ (\cite{journe}), and denote it
$\Vtilde$. The associated semigroups of the dual cocycle are related
to those of $V$ as follows:
\begin{equation} \label{dual semigroups}
\wt{Q}^{c,d}_t = (Q^{d,c}_t)^*.
\end{equation}

\subsection{Continuity} 

The above proposition makes no continuity demands on the time variable
of $V$ --- indeed it does not even require the weak measurability
condition imposed on processes. However the decomposition of a cocycle
into its associated semigroups does provide a useful handle on the
continuity of a cocycle.

\begin{lemma} \label{cts coc}
Let $V$ be a left contraction cocycle on $\init$ and let $\{\Qcd: c,d
\in \noise\}$ be its associated semigroups. Then, the following are
equivalent\tu{:}
\begin{rlist}
\item \label{coc iv}
$V$ is strongly continuous.
\item \label{coc i}
$V$ is weakly continuous at $0$.
\item \label{coc iii}
For all $c,d \in \noise$, $\Qcd$ is strongly continuous.
\item \label{coc ii}
For some $a,b \in \noise$, $\Qab$ is weakly continuous at $0$.
\end{rlist}
\end{lemma}

\begin{proof}
Let $c,d \in \noise$ and suppose that $V$ is weakly continuous at $0$.
By adaptedness
\[
\ip{u}{\Qcd_t v} = \ip{\w{c\tT}}{\w{d\tT}}^{-1} \ip{u \w{c\toT}}{V_t v
\w{d\toT}},
\]
for $0 \leqslant t \leqslant T$ and $u,v \in \init$, and so the
contraction semigroup $\Qcd$ is weakly continuous at $0$ and thus also
strongly continuous, by standard semigroup theory (\cite{Davies},
Proposition~1.23). Thus~\eqref{coc i} implies~\eqref{coc iii}. Suppose
now that $\Qab$ is weakly continuous at $0$ (and thus strongly
continuous) for some $a,b \in \noise$, and let $t \geqslant r$ and
$\xi = v\w{f}$ for $v \in \init$ and $f \in L^2 (\Rplus; \noise)$.
Then by the cocycle relation, and contractivity and adaptedness of
$V$,
\begin{align*}
\norm{V_t \xi - V_r\xi}^2
&\leqslant \norm{\sigma_r (V_{t-r}) \xi - \xi}^2 \\
&\leqslant 2 \re \bip{\xi}{(I-\sigma_r (V_{t-r})) \xi} \\
&= 2 \re \bip{S^*_r v \w{f_{[r,t[}}}{(I -V_{t-r}) S^*_r v
\w{f_{[r,t[}}}.
\end{align*}
For any $e \in \noise$ let $\varphi^e_{r,t} = \w{f_{[r,t[}} -
\w{e_{[r,t[}}$, so $S^*_r v \w{f_{[r,t[}} = S^*_r v \varphi^e_{r,t} +
v \w{e_{[0,t-r[}}$, and thus the right-hand side of the above is no
larger than
\[
4 \norm{v}^2 \Bigl\{ \norm{\varphi^a_{r,t}} \norm{\varphi^b_{r,t}} +
\norm{\varphi^a_{r,t}} + \norm{\varphi^b_{r,t}} \Bigr\} + 2 \Bigl|
\bip{v \w{a_{[0,t-r[}}}{(I-V_{t-r}) v \w{b_{[0,t-r[}}} \Bigr|.
\]
The first term converges to $0$ as $t-r \Tends 0$, and the second term
equals
\[
2 \Bigl| \bip{v}{(I -\Qab_{t-r}) v} + \norm{v}^2 \Bigl\{ \exp (r-t)
\chi(a,b) -1 \Bigr\} \Bigr|,
\]
hence $\norm{ (V_t-V_r) \xi} \Tends 0$ as $t-r \Tends 0$ by the
assumption on $\Qab$. Therefore, since the collection of such vectors
$\xi$ is total, the uniform boundedness of $V$ implies that it is
strongly continuous. Thus~\eqref{coc ii} implies~\eqref{coc iv}.
Since the implications~\eqref{coc iv} $\Implies$ \eqref{coc i}
and~\eqref{coc iii} $\Implies$ \eqref{coc ii} are trivial the proof is
complete.
\end{proof}

Thus, as for semigroups, strong continuity for a left contraction
cocycle is equivalent to weak continuity at $0$, and also to any of
its associated semigroups --- in particular its \dfn{Markov semigroup}
$Q^{0,0}$ --- being a $C_0$-semigroup.

\begin{rem}
By the strong continuity of $t \mapsto R_t$, Lemma~\ref{cts coc} is
equally true for right contraction cocycles.
\end{rem}

Suppose that $V$ is a strongly continuous left contraction cocycle on
$\init$. Then each $\VQcd$ is a $C_0$-contraction semigroup by
Lemma~\ref{cts coc} and so has a generator $\Gcd^V$. For immediate
purposes it is convenient to also work with the $C_0$-semigroups
defined by
\[
\VPcd_t := E^{\e{c\tot}} V_t\ E_{\e{d\tot}} = e^{t (\norm{c}^2
+\norm{d}^2)/2} \ \VQcd_t, \quad c,d \in \noise,
\]
and their generators, which we denote $\HVcd$. The generators are
related by
\[
\HVcd - \ip{c}{d} = \Gcd^V + \chi (c,d),
\]
with equality of domains, where $\chi$ is the function defined
in~\eqref{chi}.

Note also that, from~\eqref{dual semigroups},
\[
\HVtildecd = (H^V_{d,c})^* \text{ and } \GVtildecd = (G^V_{d,c})^*.
\]

\subsection{Operators associated with a cocycle}

In this section the generators of semigroups associated with a
strongly continuous left contraction cocycle are compared. This will
lead to natural sufficient domain conditions for such a cocycle to be
governed by a QSDE. First note two consequences of the
estimates~\eqref{basic} and~\eqref{refined}. For locally bounded
functions $f$ and $g$ in $L^2 (\Rplus; \noise)$,
\begin{align} \label{basic Holder}
& \norm{\e{f_{[s,t[}}- \e{g_{[s,t[}}} = O \bigl(\sqrt{t-s}\bigr)
\text{ and } \notag \\
&\norm{\e{f_{[s,t[}} - \e{g_{[s,t[}} - (f-g)_{[s,t[}} =
O (t-s),
\end{align}
as $(t-s) \Tends 0$ with $[s,t[$ in some finite interval. In
particular, for $a,c \in \noise$,
\begin{align}
&\norm{\e{c_{[s,t[}} - \e{a_{[s,t[}}} = O \bigl( \sqrt{t-s} \bigr)
\text{ and } \label{basic const} \\
&\norm{\e{c_{[s,t[}} - \e{a_{[s,t[}} - (c-a)_{[s,t[}} = O(t-s).
\label{refined const}
\end{align}
The former refines to
\begin{equation} \label{basic divided}
(t-s)^{-1/2} \norm{\e{c_{[s,t[}} - \e{a_{[s,t[}}} \leqslant \norm{c-a}
+O(t-s).
\end{equation}

Viewing $\init \ot \noise \ot L^2 (\Rplus) = \init \ot L^2 (\Rplus;
\noise)$ as a subspace of $\init \ot \Fock$, define two families of
operators associated with a cocycle $V$:
\begin{align*}
T^V_d (t) :&= t^{-1} E^{1\tot} V_t E_{\e{d\tot}}, \text{ and } \\
C^V(t) :&= t^{-1} E^{1\tot} V_t E_{1\tot},
\end{align*}
for $d \in \noise$ and $t >0$. Thus $T^V_d(t) \in B(\init; \init \ot
\noise)$, and $C^V(t) \in B(\init \ot \noise)$ is a contraction. Since
its associated semigroups satisfy
\begin{multline*}
\bip{u}{(\VPcd_t -\VPad_t) v} -t \ip{u(c-a)}{T^V_d (t) v} \\
= \bip{u \bigl(\e{c\tot} -\e{a\tot} -(c-a)\tot \bigr)}{V_t v
\e{d\tot}- v\e{0}}
\end{multline*}
and $V_t v \e{d\tot} \Tends v \e{0}$ as $t \Tends 0^+$, the
estimate~\eqref{refined const} implies that
\begin{equation} \label{IP T}
t \ip{u(c-a)}{T^V_d (t)v} = \ip{u}{(\VPcd_t - \VPad_t)v} + o(t)
\end{equation}
as $t \Tends 0$; similarly
\begin{equation} \label{IP C}
t \bip{u(c-a)}{C^V (t) v(d-b)} = \bip{u}{(\VPcd_t -\VPcb_t -\VPad_t
+\VPab_t) v} +O(t^{3/2}).
\end{equation}
Now define operators $T^V_d$ and $C^V$ by
\[
T^V_d v = \underset{t \Tends 0^+}{\wlim} \, T^V_d (t) v \ 
\text{ and } \ C^V \xi = \underset{t \Tends 0^+}{\wlim} \, C^V(t) 
\xi
\]
with domains equal to the subspaces on which weak convergence holds.
Thus $\Dom C^V$ is a closed subspace of $\init \ot \noise$ on which
$C^V$ is a contraction. We shall see that each $T^V_d$ is densely
defined and obtain sufficient conditions for $C^V$ to be defined on
all of $\init \ot \noise$.

For a subset $S$ of $\noise$ and element $d$ of $\noise$, define
\begin{equation} \label{DVS}
\Dvd := \Dom G^V_{d,d}, \ \text{ and } \ \D^{V,S} := \bigcap_{c \in S}
\D^{V,c}.
\end{equation}
Note that $\Dom H^V_{c,d} = \Dom G^V_{c,d}$ for all $c,d \in \noise$.

In the next result we shall spare the reader a panoply of symbols by
dropping the $V$ and $\Vtilde$ superscripts, writing $\wt{T}_c$ for
$T^\Vtilde_c$ and so forth.

\begin{propn} \label{generators compared}
Let $V$ be a strongly continuous left contraction cocycle on $\init$
with noise dimension space $\noise$, let $c,d \in \noise$.
\begin{alist}
\item
For each $v \in \Dd$ and $f \in (L^2 \cap L^\infty_\loc) (\Rplus;
\noise)$, the map $t \mapsto V_t v \e{f}$ is \tu{(}locally\tu{)}
H\"{o}lder continuous with exponent $\tfrac{1}{2}$.
\item
For each $a \in \noise$, $\Dom H_{a,d} = \Dd$. Moreover, for any dense
subspace $\D$ of $\Dd$, $\D$ is a core for $H_{c,d}$ if and only if it
is a core for $H_{d,d}$.
\item
For each $e \in \noise$,
\begin{equation} \label{EeTd}
E^e T_d \supset H_{c+e,d} - H_{c,d}.
\end{equation}
In particular, the map $c \mapsto H_{c,d}$ is complex-conjugate affine
linear.
\item
For all $b,e \in \noise$,
\begin{equation}
\label{HcdHcb}
H_{c,d} -H_{c,b} \subset (\wt{T}_c)^* E_{d-b} \quad \text{and} \quad
H_{c,d} +(\wt{T}_c)^* E_e \subset H_{c,d +e},
\end{equation}
and the maps $d \mapsto H_{c,d}$ and $d \mapsto T_d$ are complex
affine linear in the sense that, for all $b \in \noise$ and $z \in
\Comp$, if $e = (1-z) b+zd$ then
\begin{equation} \label{H and T incl}
H_{c,e} \supset (1-z) H_{c,b} + zH_{c,d} \quad \text{and} \quad T_e
\supset (1-z) T_b+zT_d.
\end{equation}
\item
For $b,e \in \noise$
\[
T_d -T_b \subset CE_{d-b} \quad \text{and} \quad T_d +CE_e \subset
T_{d+e},
\]
in particular the operator $(T_d -T_b)$ is bounded on its domain. If
$\Dom CE_d = \init $ then $(CE_d)^* \supset E^d \wt{C}$. Also $\Dom C=
\init \ot \noise$ if and only if $\Dom \wt{C} = \init \ot \noise$, in
which case $\wt{C} = C^*$ and $\Dom T_d$ is independent of $d$.
\item
For each $b \in \noise$, $T_b$ is $H_{c,d}$-bounded, with relative
bound $0$, on $\Dd \cap \Dom T_b \supset \D^{\{b,d\}}$, in the
notation~\eqref{DVS}.
\end{alist}
\end{propn}

\begin{proof}
First note that for $v \in \Dd$, and $\lambda, t >0$,
\begin{align}
t^{-1} \norm{(V_t -I) v\w{d\tot}}^2 \notag
&\leqslant 2t^{-1} \re \bip{v}{(I -Q^{d,d}_t) v} \notag \\
&= -2t^{-1} \re \int^t_0 \ip{v}{Q^{d,d}_s G_{d,d} v} \, ds \notag \\
&\leqslant 2\norm{v} \norm{G_{d,d} v} \notag \\
&\leqslant 2 \bigl( \norm{v} \norm{H_{d,d} v} + \norm{d}^2 \norm{v}^2
\bigr) \notag \\
&\leqslant \Bigl( \lambda \norm{H_{d,d}v} + \bigl( \lambda^{-1} +
\sqrt{2} \norm{d} \bigr) \norm{v} \Bigr)^2, \notag \\
\intertext{so}
t^{-1/2} \norm{(V_t-I) v\e{d\tot}} &\leqslant \lambda \norm{H_{d,d}v}
+ \mu(\lambda) \norm{v} + O(t) \label{Vt-I}
\end{align}
as $t \Tends 0$, where $\mu(\lambda) = \lambda ^{-1} + \sqrt{2}
\norm{d}$. In particular, since $E^{1\tot}E_{\e{0}} =0$,
using~\eqref{basic divided}
\begin{align} \label{Tdv bdd}
\norm{T_d(t) v}
&\leqslant t^{-1/2} \norm{V_t v \e{d\tot} -v\e{0}} \notag \\
&\leqslant \lambda \norm{H_{d,d} v} + \mu'(\lambda) \norm{v} + O(t)
\end{align}
where $\mu'(\lambda) = \lambda^{-1} + (1+\sqrt{2}) \norm{d}$.
From~\eqref{IP T} therefore
\begin{align} \label{P-P}
\norm{(P^{c,d}_t - P^{a,d}_t) v}
&\leqslant t \norm{c-a} \norm{T_d(t)v} + o(t) \notag \\
&\leqslant t \norm{c-a} \bigl( \lambda \norm{H_{d,d}v} + \mu'(\lambda)
\norm{v} \bigr) + o(t).
\end{align}

(a) Let $v \in \Dd, f \in (L^2 \cap L^\infty_\loc) (\Rplus; \noise)$
and $T \geqslant t \geqslant s \geqslant 0$. Then, from the
estimate~\eqref{basic Holder} the function $f_{s,t} := f+
(d-f)_{[s,t[}$ satisfies
\[
\norm{\e{f} - \e{f_{s,t}}} = O(\sqrt{t-s}),
\]
as $(t-s) \Tends 0$. Using the cocycle and adaptedness properties of
$V$,
\begin{align*}
\norm{(V_t -V_s) v\e{f_{s,t}}} &\leqslant \norm{\sigma_s (V_{t-s} -I)
v \e{f_{s,t}}} \\
&= \norm{\e{f_{[s,t[^c}}} \norm{\sigma_s (V_{t-s} -I) v\e{d_{[s,t[}}}
\\
&\leqslant \norm{\e{f}} \norm{(V_{t-s} -I) v\e{d_{[0,t-s[}}}.
\end{align*}
Therefore, by~\eqref{Vt-I}, if follows that $\norm{(V_t-V_s) v\e{f}} =
O(\sqrt{t-s})$, and (a) follows.

\medskip
(b) By~\eqref{P-P}
\[
\limsup_{t \Tends 0^+} t^{-1} \norm{(P^{c,d}_t - P^{d,d}_t)v} < \infty
\text{ for } v \in \Dd.
\]
By standard semigroup theory (\cite{Davies}, Corollary~1.39) it
therefore follows that $\Dom H_{c,d} \supset \Dd$, and~\eqref{P-P}
gives
\[
\norm{(H_{c,d} - H_{d,d}) v} \leqslant \norm{c-d} \bigl( \lambda
\norm{H_{d,d} v} + \mu'(\lambda) \norm{v} \bigr), \quad v \in \Dd.
\]
We may therefore apply Gustafson's Theorem (Theorem~\ref{Gustafson})
with $A= H_{c,d} |_{\Dd}$, $B=H_{d,d}$ and an appropriately chosen
$\lambda$ to conclude that
\[
\Dom \ol{A} = \Dom \ol{B} = \Dom B = \Dom A
\]
and that $A$ itself generates a $C_0$-contraction semigroup. But
$C_0$-contraction semigroup generators are maximal dissipative
(\cite[Theorem 6.4]{Davies} or \cite[page 241]{ReS}) so the inclusion
$H_{c,d} \supset A$ is an equality --- in other words $\Dom H_{c,d} =
\Dd$. A further application of Gustafson's Theorem now shows that, for
a dense subspace $\D$ of $\Dd$, $\D$ is a core for $H_{c,d}$ if and
only if it is a core for $H_{d,d}$.

\medskip
(c) Let $v \in \Dd$. First note that, by~\eqref{Tdv bdd}, $T_d(t) v$
is locally bounded in $t$ in a neighbourhood of $0$. Thus, in view
of~(b) and~\eqref{IP T}, $v \in \Dom T_d$ and $E^e T_d v = (H_{c+e,d}
- H_{c,d})v$ for any $e \in \noise$. Thus~\eqref{EeTd} holds.

\medskip
(d) Let $b\in \noise, v \in \D^{\{b,d\}}$ and $u \in \Dom \wt{T}_c$.
Then, for $t >0$,
\[
\bip{v (d-b)}{\wt{T}_c(t) u} = t^{-1} \bip{v}{(\wt{P}^{d,c}_t -
\wt{P}^{b,c}_t) u} + o(1) = t^{-1} \bip{(P^{c,d}_t - P^{c,b}_t)v}{u} +
o(1).
\]
It follows that $v (d-b) \in \Dom (\wt{T}_c)^*$ and $(\wt{T}_c)^*
E_{d-b} v =(H_{c,d} - H_{c,b}) v$. This proves the first inclusion
in~\eqref{HcdHcb}; \eqref{EeTd} applied to $\Vtilde$ gives
$\wt{H}_{d,c} + E^e \wt{T}_c = \wt{H}_{d+e,c}$ which yields the
second:
\[
H_{c,d+e} = (\wt{H}_{d,c} + E^e \wt{T}_c)^* \supset H_{c,d} +
\wt{T}^*_c E_e.
\]
By (c), setting $e= (1-z) b+zd$,
\[
(1-z) H_{c,b} +z H_{c,d} \subset \bigl( (1 -\ol{z}) \wt{H}_{b,c}
+\ol{z} \wt{H}_{d,c} \bigr)^* = (\wt{H}_{e,c})^* = H_{c,e}.
\]
This gives the first of the inclusions~\eqref{H and T incl}; the
second follows from the observation
\begin{align*}
\norm{(1-z)T_b (t) +z T_d (t) -T_e (t)}
&\leqslant t^{-1/2} \norm{(1-z) \e{b\tot} +z\e{d\tot} - \e{e\tot}} \\
&= O (t^{1/2}).
\end{align*}

\medskip
(e) The first two inclusions follow from the observation
\begin{align*}
\norm{T_d(t) -T_b (t) -C(t)E_{d-b}}
&\leqslant t^{-1/2} \norm{\e{d\tot} -\e{b\tot}-(d-b)\tot} \\
&=O(t^{1/2}),
\end{align*}
by~\eqref{refined const}, and the rest follows from the fact that
$\wt{C}(t) = C(t)^*$ for each $t >0$.

\medskip
(f) Let $v \in \Dd$ and $\lambda >0$. From~\eqref{P-P}
\[
\norm{(H_{c,d} - H_{a,d}) v} \leqslant \norm{c-a} \bigl( \lambda
\norm{H_{d,d} v} + \mu' (\lambda) \norm{v} \bigr)
\]
Taking $a=d$ it follows that
\[
\bigl(1 -\norm{c-d} \lambda \bigr) \norm{H_{d,d} v} \leqslant
\norm{H_{c,d} v} + \norm{c-d} \mu'(\lambda) \norm{v}.
\]
But from~\eqref{Tdv bdd} it follows that
\[
\norm{T_d v} \leqslant \lambda \norm{H_{d,d} v} + \mu'(\lambda)
\norm{v},
\]
therefore $T_d$ is $H_{c,d}$-bounded with relative bound $0$. Since
$(T_d -T_b)$ is bounded on its domain~(f) follows. This completes the
proof.
\end{proof}

To a strongly continuous left contraction cocycle $V$ on $\init$, with
noise dimension space $\noise$, we may therefore associate an operator
on $\init \oplus (\init \ot \noise)$ by
\[
F^V := \begin{bmatrix} Z^V & M^V \\ L^V & C^V-I \end{bmatrix}
\]
where $Z^V=H^V_{0,0}$, $L^V =T^V_0 |_{\Dvo}$ and $M^V =(T^\Vtilde_0
|_{\D^{\Vtilde, 0}})^*$. Thus $Z^V$ is a $C_0$-contraction semigroup
generator, $L^V$ has the same dense domain as $Z^V$, $M^V$ is closed
and $C^V$ is a contraction operator. If $F^V$ is densely defined then
\begin{equation} \label{FV*}
(F^V)^* \supset F^\Vtilde.
\end{equation}

\begin{cor} \label{C1.4}
For all $c,d \in \noise$ and $S \subset \noise$
\[
\Dom G^V_{c,d} = \Dvd \ \text{ and } \ \D^{V, \Aff S} = \D^{V,S}
\]
where $\Aff S$ denotes the complex affine span of $S$. Moreover,
\begin{equation} \label{1.19b}
\Dom F^V \supset \Dvo \op \bigl( \D^{V, \Aff S} \uot \iD \bigr),
\end{equation}
where $\iD = \Lin (S-S)$.
\end{cor}

\begin{proof}
For convenience we drop the superscripts $V$ and $\Vtilde$ as in the
proposition. The semigroup generators $\Gcd$ and $H_{c,d}$ have
the same domains, so the first equality follows from part~(b) of the
proposition. For the second equality, if $e \in \Aff S$ then
$\D^e \supset \D^S$ by the first inclusion in~\eqref{H and T incl}.
But this implies that
\[
\D^S \subset \bigcap_{e \in \Aff S}\, \D^e = \D^{\Aff S} \subset
\D^S.
\]
For~\eqref{1.19b} note that if $b,d \in \noise$ then
\[
ME_{d-b} = (E^{d-b} \wt{T}_0 |_{\D^{0}})^* = (\wt{H}_{d,0} -
\wt{H}_{b,0})^* \supset H_{0,d} - H_{0,b},
\]
applying part~(c) to the dual cocycle $\Vtilde$. Also $CE_{d-b}
\supset T_d - T_b$, thus $\D^{\{b,d\}}$ is a subspace of both $\Dom
ME_{d-b}$ and $\Dom CE_{d-b}$. Therefore if $e \in S-S$ and $v \in
\D^S$ then $ve \in \Dom M \cap \Dom C$. The result follows since $\Dom
L = \Dom Z=\D^0$.
\end{proof}

For a cocycle $V$ and subspace $\iD$ of $\noise$, Corollary~\ref{C1.4}
permits the following definition:
\begin{equation} \label{FVT}
F^{V,\iD} := F^V |_{\D_0 \op (\D \uot \iD)}
\end{equation}
where $\D_0 = \Dvo $ and $\D = \DvD$. 
Note that for any subset $\iT$ of $\noise$ containing $0$
\begin{equation} \label{DVT=DVD}
\D^{V,\iT} = \D^{V,\iD}, \ \text{ where } \iD = \Lin \iT.
\end{equation}

From the corollary we see that $\D^{V, \{0,d\}} \subset \Dom
E^{\wh{c}} F^V E_{\wh{d}}$ for all $c,d \in \noise$, and by
parts~(b)--(e) of the proposition,
\begin{equation} \label{FV,H,G}
E^{\wh{c}} F^V E_{\wh{d}} = H^V_{c,d} - \ip{c}{d} = G^V_{c,d} +\chi
(c,d) \text{ on } \D^{V,\{d,0\}}.
\end{equation}

For \dfn{Markov-regular} cocycles, that is cocycles whose Markov
semigroup $Q^{0,0}$ is norm-continuous, the situation is much simpler.

\begin{cor} \label{Mreg cocs}
Let $V$ be a strongly continuous left contractive cocycle on $\init$
with noise dimension space $\noise$ and suppose that one of its
associated semigroups $\Qcd$ is norm continuous. Then all of its
associated semigroups are norm continuous and $F^V \in B(\init \ot
\khat)$.
\end{cor}

\begin{proof}
That all or none of the associated semigroups are norm continuous
follows since
\[
\norm{P^{a,b}_t - P^{c,d}_t} = O(\sqrt{t}) \ \text{ for } a,b,c,d \in
\noise,
\]
by~\eqref{basic const}. So if it \emph{is} the case that all the
semigroups are norm continuous then $H^V_{c,d} \in B(\init)$ for all
$c,d \in \noise$. In particular $Z^V = H^V_{0,0} \in B(\init)$, so
that $\D^{V,0} = \init$, and hence $L^V = T^V_0$ which is bounded by
part~(f) of the proposition.

Since the semigroups associated to the dual cocycle are the adjoints
of those associated with the cocycle, and must also be norm
continuous, $M^V = T^{\Vtilde*}_0 \in B(\init \ot \noise; \init)$, and
from~\eqref{IP C} it follows that the contraction $C^V$ is densely
defined, hence $C^V\in B(\init \ot \noise)$.
\end{proof}

\subsection{Cocycle characterisation through semigroups}

In Proposition~\ref{Semi decomp} we used the family of maps
$\{\VQcd_t: c,d \in \iT, t \geqslant 0\}$ defined in terms of a
\ti{given} process $V$ to determine whether or not it is a left
cocycle. The following result turns this around.

\begin{thm}[\cite{spawn}] \label{recon thm}
Let $\Q_\iT = \{\Qcd: c,d \in \iT\}$ be a family of semigroups on
$\init$ indexed by a total subset $\iT$ of $\noise$ which contains
$0$. Then the following are equivalent\tu{:}
\begin{rlist}
\item
There is a left contraction cocycle $V$ on $\init$ whose associated
family of semigroups includes $\Q_\iT$.
\item
For all $n \in \Nat$, $Y \in \Mat_n (|\init\rangle) = B(\Comp^n;
\init^n)$, and positive invertible matrices $A, B \in \Mat_n (\Comp)$,
if $\norm{A^{-1/2} Y B^{-1/2}} \leqslant 1$ then
\begin{equation} \label{recon ineqs}
\norm{(A \schur \bvarpi^\bfc_t)^{-1/2} (Q^\bfc_t \schur Y) (B
\schur \bvarpi^\bfc_t)^{-1/2}} \leqslant 1,
\end{equation}
for all $\bfc \in \iT^n, t \geqslant 0$.
\end{rlist}
\end{thm}

This requires some explanation of terms: $|\init\rangle := B(\Comp;
\init)$, the column operator space determined by $\init$ (\cite{EfR,
Pisier}); given $c,d \in \noise$, $\varpi^{c,d}_t =
\ip{\w{c\tot}}{\w{d\tot}}= \exp -t \chi(c,d)$, and given $\bfc \in
\noise^n$, $\bvarpi^\bfc_t := [\varpi^{c_i,c_j}_t] \in \Mat_n
(\Comp) =B(\Comp^n)$, $Q^\bfc_t := [Q^{c_i,c_j}_t] \in \Mat_n
(B(\init)) = B(\init^n)$, the symbol $\schur$ denotes the Schur
product of matrices, so in particular if $T = [|u^i_j\rangle]$ then
$Q^\bfc_t \schur T = [|Q^{c_i,c_j}_t u^i_j\rangle] \in B(\Comp^n;
\init^n)$; finally, the first matrix within each of the norms is
thought of as having entries of the form $\nu I_\init$ for $\nu \in
\Comp$, thus both norms are those of $B(\Comp^n; \init^n)$.

What this result tells us is that if we can find a family of
semigroups $\Q_\iT$ on $\init$, indexed by a total subset $\iT$ of
$\noise$ containing $0$, which satisfies~\eqref{recon ineqs} then
there is an associated cocycle $V$ on $\init$. This condition, on a
putative family of semigroups $\Q_\iT$, looks hard to verify. However
the strength of the result lies in the fact that it is manifestly
stable under pointwise limits.

\begin{thm} \label{semi approx}
Let $\Q_\iT = \{\Qcd: c,d \in \iT \}$ be a family of semigroups on
$\init$, indexed by a total subset $\iT$ of $\noise$ which includes
$0$. Suppose that there is a sequence $(V^{(n)})_{n \geqslant 1}$ of
left contraction cocycles on $\init$ whose associated semigroups
satisfy
\[
\Qncd_t \Tends \Qcd_t \quad \text{ pointwise on } \init,
\]
for all $c,d \in \iT$ and $t >0$. Then there is a unique left
contraction cocycle $V$ on $\init$ whose associated semigroups include
the family $\Q_\iT$. Moreover $V^{(n)}_t \Tends V_t$ in the weak
operator topology for each $t$. 
\end{thm}

\begin{proof}
The existence of $V$ is immediate from Theorem~\ref{recon thm},
uniqueness follows from the totality of $\iT$, and the convergence
$V^{(n)} \Tends V$ is a consequence of~\eqref{semi decomp} and
contractivity.
\end{proof}

\section{Quantum Stochastic Differential Equations}

Let $\khat := \Comp \op \noise$ and, for any subspace $\iD$ of
$\noise$, let $\iDhat = \Comp\op \iD = \Lin\{\wh{d}: d \in \iD\}$
where $\wh{d} := \binom{1}{d}$. Also let $e_0 = \binom{1}{0} \in
\khat$ and define
\[
\Delta = I_\init \ot P_\noise
\]
where $P_\noise \in B (\khat)$ is the orthogonal projection with range
$\{e_0\}^\perp = 0 \oplus \noise$.

Consider now the \dfn{form QSDE}
\begin{equation} \label{form LHP}
dV_t = V_t \, \Form \, d\Lambda_t, \quad V_0 =I,
\end{equation}
and the \dfn{operator QSDE}
\begin{equation} \label{op LHP}
dV_t = \wh{V}_t (F \uot I_\Fock) \, d\Lambda_t,\quad V_0 = I,
\end{equation}
for (bounded operator-valued) processes $V$ on $\init$, for which we
need dense subspaces $\D_0 \supset \D$ of $\init$, and total subsets
$\iT^\dagger$ and $\iT$ of $\noise$ containing $0$. Set $\iD = \Lin
\iT $ and $\iD^\dagger= \Lin \iT^\dagger$.

In the first case, $\Form$ is an operator-valued map defined on
$\wh{\iD^\dagger} \times \iDhat$ of the form
\[
\left( \binom{z}{c}, \binom{w}{d} \right) \mapsto \begin{bmatrix}
\ol{z} & 1 \end{bmatrix} \begin{bmatrix} K & M_d \\ L^c & N^c_d
\end{bmatrix} \begin{bmatrix} w \\ 1 \end{bmatrix} = \ol{z}w K +\ol{z}
M_d +w L^c +N^c_d
\]
where $K \in \Op (\D_0)$, $c \mapsto L^c$ is conjugate linear
$\iD^\dagger \To \Op (\D_0)$, $d \mapsto M_d$ is linear $\iD \To \Op
(\D)$ and $(c,d) \mapsto N^c_d$ is sesquilinear $\iD^\dagger \times
\iD \To \Op (\D)$ (thus $\Form (\binom{z}{c}, \binom{w}{d}) \in \Op
(\D)$ in general, and is $\Op (\D_0)$-valued if $d=0$), and $V$ is a
$\iT^\dagger$-$\iT$-\dfn{solution} of~\eqref{form LHP} on $\D_0 \ot
\e{0} +\D \uot \Exp_\iT$ if, in the notation $\wh{g}(s):=\wh{g(s)}$,
\begin{equation} \label{form soln}
\ip{u\e{f}}{(V_t -I) v\e{g}} = \int^t_0 \ip{u\e{f}} {V_s \Form
(\wh{f}(s), \wh{g}(s)) v\e{g}} \, ds
\end{equation}
for all $u \in \init$, $f \in \step_{\iT^\dagger}$, $(v,g) \in \bigl(
\D_0 \times \{0\} \bigr) \cup \bigl( \D \times \step_\iT \bigr)$ and
$t \geqslant 0$. In particular $V$ is weakly continuous in an obvious
sense.

In the second case $F \in \Op \bigl(\D_0 \op (\D \uot \iD)\bigr)$,
$\wh{V}_t$ stands for the operator on $\init \ot \khat \ot \Fock$
obtained from $V_t \ot I_{\khat}$ by tensor flipping, and there are
two basic kinds of solution: $V$ is a $\iT^\dagger$-\dfn{weak
solution} of~\eqref{op LHP} on $\D_0 \ot \e{0} +\D \uot \Exp_\iT$
if~\eqref{form soln} holds for the \dfn{component map} of $F$, defined
by
\[
\Form (\xi, \eta) = E^\xi FE_\eta, \quad \xi \in \wh{\iD^\dagger},
\eta \in \iDhat.
\]
In other words, setting $\zeta (s) = v e_0 \e{0} +w \wh{g} (t) \e{g}$, 
\begin{multline*}
\bip{u \e{f}}{(V_t -I) \{v\e{0} +w\e{g}\}} \\
\begin{aligned}
&= \int^t_0 \bip{u\e{f}} {V_s E^{\wh{f}(s)} F \uot I_\Fock \zeta (s)}
\, ds \\
&= \int^t_0 \bip{u\e{f}}{V_s \{ (K + L^{f(s)}) v\e{0} + (K+L^{f(s)}
+M_{g(s)} +N^{f(s)}_{g(s)}) w\e{g} \}} \, ds,
\end{aligned}
\end{multline*}
where $F = \begin{sbmatrix} K & M \\ L & N \end{sbmatrix}$ in block
matrix form, and $L^c = E^c L$ etc.

$V$ is a \dfn{strong solution} of~\eqref{op LHP} on the same domain if
the map
\[
t \mapsto \norm{\wh{V}_t \Delta F \uot I_\Fock \zeta (t)}^2 +
\norm{\wh{V}_t \Delta^\perp F \uot I_\Fock \zeta(t)}
\]
is locally integrable for each $v \in \D_0$, $w \in \D$, and $g \in
\step_\iT$, $V$ is strongly measurable, (and hence is stochastically
integrable), and if it satisfies the quantum stochastic integral
equation
\[
V_t = I + \int^t_0 \wh{V}_s (F \uot I_\Fock) \, d\Lambda_s.
\]
In particular $V$ is strongly continuous on its domain (\cite{L}),
hence on all of $\init \ot \Fock$ by contractivity, and the First
Fundamental Formula of quantum stochastic calculus implies that $V$ is
a $\noise$-weak solution on $\D_0 \ot \e{0} +\D \uot \Exp_\iT$.

We recall the basic implication for $F$ of contractivity of a strong
solution of~\eqref{op LHP}, and include its short proof for the
convenience of the reader (cf.~\cite{RHP}).

\begin{propn}[\cite{FF8, MoFest}] \label{nec ineq}
Let $F \in \Op \bigl(\D_0 \oplus (\D \uot \iD) \bigr)$ where $\D_0
\supset \D$ are dense subspaces of $\init$ and $\iD = \Lin\iT$ for a
total subset $\iT$ of $\noise$ containing $0$, and suppose
that~\eqref{op LHP} has a strong contractive solution on $\D_0 \ot
\e{0} + \D \uot \Exp_\iT$. Then $F$ satisfies the form inequality
\begin{equation} \label{form ineq}
2 \re \ip{\xi}{F\xi} +\norm{\Delta F \xi}^2 \leqslant 0,
\end{equation}
with equality if the solution is isometric.
\end{propn}

\begin{proof}
Let $\xi \in \D_0 \oplus \bigl(\D \uot \iD\bigr)$. Then $\xi$ is
expressible in the form $u_0 e_0 + \sum^n_{i=1} u_i \wh{c}_i$ for some
$u_0 \in\D_0, n \in \Nat, \, u_1, \ldots, u_n \in \D$ and $c_1,
\ldots, c_n \in \iT$. Let $\zeta = u_0 \e{0} + \sum^n_{i=1} u_i
\e{f_i}$ and $\zeta (s) = u_0 e_0 \e{0} + \sum^n_{i=1} u_i \wh{f}_i
(s) \e{f_i}$ where $f_i = {c_i}\toT$ for $i=1, \ldots, n$ and some $T
> 0$. Then, by the Second Fundamental Formula of quantum stochastic
calculus,
\begin{align*}
0 &\geqslant t^{-1} \bigl( \norm{V_t \zeta}^2 - \norm{\zeta}^2 \bigr)
\\
&= t^{-1} \int^t_0 \bigl\{ 2 \re \ip{\wh{V}_s \zeta (s)}{\wh{V}_s (F
\uot I_\Fock) \zeta (s)} +\norm{\wh{V}_s (\Delta F \uot I_\Fock) \zeta
(s)}^2 \bigr\} \, ds,
\end{align*}
with equality if $V$ is isometric. Using the continuity of the
integrand at the origin, letting $t \Tends 0$ and then letting $T
\Tends 0$ now gives the result.
\end{proof}

\begin{rems}
(i) If it is assumed further that all the $\iT$-components
$(E^{\wh{c}} FE_{\wh{d}}: c,d \in \iT)$ are bounded, then~\eqref{op
LHP} may be solved by Picard iteration and Mohari and Fagnola showed
that~\eqref{form ineq} is also \emph{sufficient} for contractivity of
the solution. In fact boundedness of $\iT$-components and
contractivity of the solution implies that $F$ itself is bounded
(\cite{mother}, Theorem~7.5) so that~\eqref{form ineq} simplifies to
the operator inequality
\begin{equation} \label{op ineq}
F +F^* +F^* \Delta F \leqslant 0.
\end{equation}
The solution is also \emph{unique} amongst $\iT^\dagger$-weak
solutions (cf.\ Theorem~\ref{uniqueness} below).

(ii) Since the integrability condition for being a strong solution is
automatically satisfied by strongly measurable contraction processes,
any strongly measurable weak solution of~\eqref{op LHP} is necessarily
a \emph{strong} solution on the same domain.
\end{rems}

We next show how the assumption of strong measurability renders form
solutions into strong operator solutions.

\begin{thm}
Let $V$ be a strongly measurable contraction process on $\init$ with
noise dimension space $\noise$, and let $\iT$ be a total subset of
$\noise$ containing $0$, let $\D_0 \supset \D$ be dense subspaces of
$\init$. Set $\iD = \Lin \iT$, and assume that $\iD$ has an
orthonormal basis. If $V$ is a $\iT$-$\iT$-solution of the form
\tu{QSDE}~\eqref{form LHP} on $\D_0 \ot \e{0} + \D \uot \Exp_\iT$, and
each map $\xi \mapsto \Form (\xi, \eta)v$ is continuous then $\Form$
is the component map of an operator $F \in \Op \bigl(\D_0 \op (\D \uot
\iD)\bigr)$ and $V$ satisfies the corresponding operator
\tu{QSDE}~\eqref{op LHP} strongly on the same domain.
\end{thm}

\begin{proof}
By the second remark above it suffices to show that $\Form$ is
necessarily the component map of an operator $F \in \Op \bigl(\D_0 \op
(\D \uot \iD) \bigr)$. For any subspace $\hil$ of $\noise$ of the form
$\Lin \iT_0$ where $\iT_0$ is a finite subset of $\iT$, define an
operator $F^\hil \in \Op \bigl(\D_0 \op (\D \uot \hil)\bigr)$ by
\[
F^\hil = \sum_{\al, \be} E_{e_\al} \Form (e_\al, e_\be) E^{e_\be}
\]
where $(e_\alpha)$ is an orthonormal basis for $\wh{\hil}$ which
includes the vector $e_0 =\binom{1}{0}$. By sesquilinearity, $F^\hil$
does not depend on the choice of basis and, for $\xi \in \wh{\hil}$
and $(v,\eta) \in \bigl(\D_0 \times \Comp e_0\bigr) \cup \bigl(\D
\times \wh{\hil}\bigr)$,
\[
E^\xi F^\hil v\eta = \Form(\xi,\eta)v.
\]
Also let $J^\hil$ be the natural isometric embedding $\init \ot \Gamma
\bigl(L^2 (\Rplus; \hil)\bigr) \To \init \ot \Fock$, obtained by
second quantisation of the inclusion map $\hil \To \noise$. Then it is
easily verified that the process $V^\hil := (J^{\hil*} V_t J^\hil)_{t
\geqslant 0}$ satisfies the QSDE $dX_t = \wh{X}_t (F^\hil \ot I) \,
d\Lambda_t, \, X_0 =I$, $\iT_0$-weakly on $\D_0 \ot \e{0} +\D \uot
\Exp_{\iT_0}$. Since $V^\hil$ is contractive and strongly measurable
it satisfies the equation strongly. Therefore, by Proposition~\ref{nec
ineq}, $F^\hil$ satisfies
\begin{equation} \label{Delta F h}
\norm{\Delta F^\hil v \eta}^2 \leqslant -2 \re \ip{v\eta}{F^\hil
v\eta} = -2\re \ip{v}{\Form (\eta, \eta)v}
\end{equation}
for $(v,\eta) \in \bigl(\D_0 \times \Comp e_0\bigr) \cup \bigl(\D
\times \wh{\hil}\bigr)$.

Now let $(d_i)_{i \in I}$ be an orthonormal basis for $\noise$ taken
from $\iD$ (with $0 \notin I$), and set $\wh{I} = \{ 0\} \cup I$, $e_0
= \binom{1}{0}$ and $e_i = \binom{0}{d_i}$ so that $(e_\alpha)_{\al
\in \wh{I}}$ is a basis for $\khat$. Then, for $(v,\eta) \in \bigl(\D
\times \Comp \binom{0}{d}\bigr) \cup \bigl(\D_0 \times \Comp e_0
\bigr)$ where $d \in \iT$ and $I_0$ a finite subset of $I$,
applying~\eqref{Delta F h} with $\hil = \Lin \bigl( \{d\} \cup \{d_i:
i \in I_0 \} \bigr)$ gives
\begin{align*}
\sum_{i \in I_0} \norm{\Form (e_i, \eta)v}^2
&= \Bigl\| \sum_{i \in I_0} E_{e_i} E^{e_i} F^\hil v \eta \Bigr\|^2 \\
&\leqslant \norm{\Delta F^\hil v \eta}^2 \leqslant -2 \re
\bip{v}{\Form (\eta, \eta) v},
\end{align*}
and so the orthogonal sum $\sum_{i \in I} E_{e_i} \Form (e_i, \eta) v$
is convergent. Thus an operator $F \in \Op \bigl(\D_0 \op (\D \uot
\iD) \bigr)$ is defined by linear extension of the prescription
\[
v\eta \mapsto \sum_{\al \in \wh{I}} E_{e_\al} \Form (e_\al, \eta) v.
\]
By the continuity assumption on $\Form$, $E^\xi F v\eta = \Form (\xi,
\eta) v$ for $\xi \in \iDhat$ and $v\eta$ as above, and it follows
that $F$ is independent of the choice of basis $(d_i)_{i \in I}$,
hence is the component map of $\Form$: $E^\xi F E_\eta = \Form (\xi,
\eta)$ for $\xi, \eta \in \iDhat$. This completes the proof.
\end{proof}

\begin{rems}
Dixmier showed that a pre-Hilbert space need not have an orthonormal
basis (\cite[V.70]{BourbakiTVS}); however the assumption on $\iD$ is
automatically satisfied if either the Hilbert space $\noise$ is
separable or the set $\iT$ contains a subset which is orthogonal and
total.

Since contraction processes satisfying the form QSDE are weakly
continuous, the strong measurability assumption is redundant when
$\init $ and $\noise$ are both separable or, in view of Lemma~\ref{cts
coc}, when the solution is a left cocycle.
\end{rems}

This connects with issues of uniqueness. 

\begin{thm} \label{uniqueness}
Let $\Form$ be a sesquilinear map $\wh{\iD^\dagger} \times \iDhat \To
\Op(\D)$, where $\D$ is a dense subspace of $\init$, $\iD^\dagger
=\Lin \iT^\dagger$ and $\iD =\Lin \iT$ for total subsets
$\iT^\dagger$ and $\iT$ of $\noise$ that contain $0$.
\begin{alist}
\item
Suppose that $\iT^\dagger = \Real \iT^\dagger$ and $\iT = \Real \iT$.
If $K := \Form (e_0,e_0)$ is a pregenerator of $C_0$-contraction
semigroup on $\init$ then the form \tu{QSDE}~\eqref{form LHP} has at
most one contractive $\iT^\dagger$-$\iT$-solution on $\D \uot
\Exp_\iT$.
\item
If the form \tu{QSDE}~\eqref{form LHP} has a unique contractive
$\iT^\dagger$-$\iT$-solution $V$ on $\D \uot \Exp_\iT$ then $V$ is a
left contraction cocycle.
\item
If the form \tu{QSDE}~\eqref{form LHP} has a
$\iT^\dagger$-$\iT$-solution $V$ on $\D \uot \Exp_\iT$ which is a left
contraction cocycle then $\D \subset \D^{V,\iD}$ and $E^{\wh{c}} F^V
E_{\wh{d}} |_\D = \Form (\wh{c}, \wh{d})$ for all $c \in \iD^\dagger,
d \in \iD$.
\end{alist}
\end{thm}

Part~(a) is Mohari's Uniqueness Theorem (\cite{Moh}). The invariance
of $\iT$ (and $\iT^\dagger$) under scalar multiplication can be 
weakened to the following:
\[
\text{for each } d \in \iT \text{ there exists } \epsilon > 0 \text{
such that } [0,\epsilon] d \subset \iT,
\]
which is sufficiently strong to allow differentiation at a crucial
stage in his argument. It is clear from Meyer's treatment (\cite[page
191]{meyer}) that the result remains valid in this generality.
Part~(b) is proved by verifying that, for each $t>0$,
\[
V^t_s = \begin{cases} V_s & \text{if } s \leqslant t, \\ V_t \si_t
(V_{s-t}) & \text{if } s >t,
\end{cases}
\]
defines a contraction process $V^t$ which also satisfies~\eqref{form
LHP}. This is easily checked by treating $\init \ot \Fock\tot$ as an
initial space and using the explicit action of shifts on exponential
vectors.

\begin{proof}[Proof of part~\tu{(c)}]
Pick $u \in \init, v \in \D, c \in \iT^\dagger$ and $d \in \iT$, and
set $f=c_{[0,1[}$ and $g=d_{[0,1[}$. Then for all $0<t<1$
\[
\ip{u \e{f}}{(V_t -I) v\e{g}} = \int^t_0 \ip{u \e{f}}{V_s \Form
(\wh{c},\wh{d}) v \e{g}} \, ds.
\]
But, for the same $t$,
\[
\ip{u}{(P_t^{c,d} -I) v} = \ip{u\e{f}}{(V_t -I) v\e{g}} e^{(t-1)
\ip{c}{d}} + (e^{t\ip{c}{d}} -1) \ip{c}{d} \ip{u}{v},
\]
and consequently 
\[
\lim_{t \Tends 0} t^{-1} \ip{u}{(P_t^{c,d} -I) v} = \ip{u}{(\Form
(\wh{c},\wh{d}) -\ip{c}{d}) v}.
\]
Since this holds for all $u \in \init$ it follows that $v \in \Dom
H^V_{c,d}$ and that $\Form (\wh{c},\wh{d}) -\ip{c}{d} \subset
H^V_{c,d}$ for all $c \in \iD^\dagger$ and $d \in \iD$ (\cite{Davies},
Theorem~1.24). Hence $\D \subset \D^{V,\iD}$ and so $E^{\wh{c}} F^V
E_{\wh{d}}|_\D = \Form (\wh{c},\wh{d})$ by~\eqref{FV,H,G}.
\end{proof}

\section{Necessary conditions for contractive solution}

In this section we explore necessary conditions on $F$ for the
existence of contractive solutions of~\eqref{op LHP}. Recall that
densely defined dissipative operators are closable with dissipative
closures.

\begin{propn} \label{ineq conseqs}
Let $F \in \Op \bigl( \D_0 \op (\D \uot \iD) \bigr)$, for dense
subspaces $\D_0 \supset \D$ of $\init$ and $\iD$ of $\noise$, have
block matrix form $\begin{sbmatrix} K & M \\ L & C-I \end{sbmatrix}$
and satisfy the form inequality~\eqref{form ineq}. For each $c \in
\noise$ and $d \in \iD$ define $\Gcd^0 := E^{\wh{c}} F E_{\wh{d}} -
\chi (c,d)$. Then
\begin{alist}
\item \label{NCI}
$C$ is a contraction and, for all $u \in \D_0$,
\begin{equation} \label{LK ineq}
\norm{Lu}^2 + 2\re \ip{u}{Ku} \leqslant 0;
\end{equation}
\item \label{diss}
$F$ and $\Gcd^0$ are dissipative\tu{;} let $\Gcd = \ol{\Gcd^0}$ and
$Z = G_{0,0} = \ol{K}$\tu{;}
\item \label{L}
$L$ is $\Gcd^0$-bounded, with relative bound $0$, on $\Dom
\Gcd^0$\tu{;}
\item \label{generator}
for each $a \in \noise$, $G^0_{a,d}$ is a relatively bounded
perturbation of $\Gcd^0$ with relative bound $0$, $\Dom G_{a,d} = \Dom
\Gcd$, and $G_{a,d}$ is a $C_0$-semigroup generator if and only if
$\Gcd$ is\tu{;}
\item \label{M*}
$\Dom M^* \supset \Dom K$ and, for all $u \in \D_0$,
\begin{equation} \label{M* ineq}
\norm{(L +CM^*)u}^2 +\norm{M^*u}^2 +2\re \ip{u}{Ku} \leqslant 0,
\end{equation}
in particular, $M^*$ is $K$-bounded with relative bound $0$\tu{;}
\item \label{Fbar}
$\ol{F}$ also satisfies~\eqref{form ineq}, moreover $\ol{F} \supset
F'$ where $F' = \begin{sbmatrix} Z & \ol{M} \\ L' & \ol{C}-I
\end{sbmatrix}$, $L'$ being the continuous extension \tu{(}in the
graph norm of $Z$\tu{)} of $L$ to $\Dom Z$.
\end{alist}
\end{propn}

\begin{proof}
For $\xi = \binom{u}{\eta} \in \D_0 \op (\D \uot \iD)$,~\eqref{form
ineq} is equivalent to
\begin{equation} \label{new form ineq}
\norm{Lu +C\eta}^2 \leqslant -2\re \ip{u}{Ku +M\eta} +\norm{\eta}^2,
\end{equation}
and setting $u=0$, respectively $\eta =0$, shows that~\eqref{NCI}
holds. Now abbreviate $E^cL$, $ME_d$ and $E^c CE_d$ to $L^c$, $M_d$
and $C^c_d$ respectively, and denote $\Dom \Gcd^0$ by $\D^0_d$, thus
\[
\Gcd^0 = K +L^c +M_d +C^c_d -\tfrac{1}{2} \norm{c}^2 -\tfrac{1}{2}
\norm{d}^2,
\]
where $\D^0_d = \D_0$ if $d=0$ and equals $\D$ otherwise.

If $\eta = ud$ where $d \in \iD$ and $u \in \D^0_d$ then~\eqref{new
form ineq} reads
\[
\norm{(L +C_d) u}^2 \leqslant -2\re \ip{u}{(K +M_d -\tfrac{1}{2}
\norm{d}^2) u},
\]
so, for $c \in \noise$ and $d \in \iD$,
\begin{align*}
-2 \re \ip{u}{\Gcd^0 u} 
&= -2 \re \bip{u}{\{(K +M_d -\tfrac{1}{2} \norm{d}^2) +(L^c +C^c_d -
\tfrac{1}{2} \norm{c}^2)\} u} \\ 
&\geqslant \norm{(L +C_d) u}^2 - 2\re \bip{uc}{(L +C_d)u} +
\norm{uc}^2 \\
&= \norm{(L +C_d -E_c)u}^2
\end{align*}
for all $u \in \D^0_d$. Thus $\Gcd^0$ is dissipative, moreover
\[
\norm{(L +C_d -E_c)u} \leqslant \sqrt{ -2 \re\ip{\lambda^{-1}
u}{\lambda \Gcd^0 u}} \leqslant \lambda \norm{\Gcd^0u} + \lambda^{-1}
\norm{u}
\]
for all $\lambda >0$ and so, since $(C_d -E_c)$ is bounded, $L$ is
$\Gcd^0$-bounded with relative bound $0$. Since $F$ is clearly
dissipative, we have established~\eqref{diss} and~\eqref{L}. Since
$G^0_{a,d} -\Gcd^0 = E^{a-c}L +C^{a-c}_d +\frac{1}{2} \bigl(
\norm{c}^2 -\norm{a}^2 \bigr)$ it also follows that $G^0_{a,d}$ is a
relatively bounded perturbation of $\Gcd^0$ with relative bound $0$,
and so~\eqref{generator} follows from Gustafson's Theorem.

Now let $u \in \D_0$. Then, from~\eqref{new form ineq}, $2 |\ip{u}{M
\eta}| \leqslant \norm{\eta}^2 -2 \re \ip{u}{Ku}$ for each $\eta \in
\D \uot D$. This implies that $u \in \Dom M^*$. Thus $\Dom M^* \supset
\D_0$ and~\eqref{new form ineq} reads $\norm{Lu +C\eta}^2 + 2\re
\ip{M^*u}{\eta} + 2\re \ip{u}{Ku} \leqslant \norm{\eta}^2$, now valid
for $u \in \D_0$ and $\eta \in \init \ot \noise$. Putting $\eta =
M^*u$ gives~\eqref{M* ineq}, in particular
\[
\norm{M^* u}^2 \leqslant 2 |\ip{u}{Ku}| \leqslant \bigl( \lambda
\norm{Ku} + \lambda ^{-1} \norm{u} \bigr)^2
\]
for $\lambda > 0$, showing that~\eqref{M*} holds. Since $F$ is densely
defined and dissipative it is closable and it is easily verified that
its closure contains $F'$, and that it inherits the
property~\eqref{form ineq} from $F$. Thus~\eqref{Fbar} holds too and
the proof is complete.
\end{proof}

\begin{rems}
(i) The form inequality~\eqref{form ineq} is therefore equivalent
to~\eqref{LK ineq} together with contractivity of $C$, $\Dom M^*
\supset \D_0$ and the following inequality holding for $u \in \D_0$
and $\eta \in \init \ot \noise$:
\[
\bigl|\ip{(M^* +C^*L) u}{\eta}\bigr|^2 \leqslant \bigl(2\re
\ip{u}{(-K)u} -\norm{Lu}^2\bigr) \bigl(\norm{\eta}^2
-\norm{C\eta}^2\bigr).
\]
If equality holds in~\eqref{form ineq} then $C$ is isometric and
$\norm{Lu}^2 + 2\re \ip{u}{Ku} = 0$ for all $u \in \D_0$, in turn, if
\emph{either} of these conditions hold then $M^* \supset -C^*L$.

(ii) In view of~(f), the proposition still holds if $K$, $L$ and $C$
are replaced by $Z$, $L'$ and $\ol{C}$ respectively, and $M$ is
replaced by the restriction of $\ol{M}$ to any dense subspace of its
domain of the form $\D'_1 \uot \iD'$.
\end{rems}

\begin{propn} \label{more ineq conseqs}
Let $F$ and $F'$ be as in \tu{Proposition~\ref{ineq conseqs}}. Suppose
that $Z$ is a generator of a $C_0$-semigroup and let $F^{(n)} :=
I^{(n)*} \ol{F} I^{(n)}$, where $I^{(n)} = \diag [J^{(n)}, \ I_{\init
\ot \noise}] \in B(\init \ot \khat)$, $J^{(n)}$ being the contraction
$(I-n^{-1} Z)^{-1}$. Then $F^{(n)}$ is bounded and its closure
satisfies the operator inequality~\eqref{op ineq}, and $F^{(n)} \Tends
F$ pointwise on $\D_0 \op (\D \uot \iD)$.
\end{propn}

\begin{proof}
Note first that $I^{(n)}$ leaves $\Dom F'$ invariant, and that $\Delta
I^{(n)} = \Delta$. Thus for $\xi \in \Dom F$, putting $\xi_n =I^{(n)}
\xi$,
\[
2 \re \ip{\xi}{F^{(n)} \xi} + \norm{\Delta F^{(n)} \xi}^2 = 2 \re
\ip{\xi_n}{\ol{F} \xi_n} + \norm{\Delta \ol{F} \xi_n}^2 \leqslant 0,
\]
by Proposition~\ref{ineq conseqs}, thus $F^{(n)}$ satisfies the form
inequality~\eqref{form ineq}. Now let $\begin{sbmatrix} K^{(n)} &
M^{(n)} \\ L^{(n)} & \ol{C}-I \end{sbmatrix}$ be the block matrix form
of $F^{(n)}$. Since $K^{(n)} = J^{(n)*} Z J^{(n)} \in B(\init)$ it
follows from Proposition~\ref{ineq conseqs} that $L^{(n)}$ and
$M^{(n)*}$ are bounded, and so $F^{(n)}$ is bounded, hence extends to
$\init \ot \khat$. Thus $F^{(n)}$ satisfies the operator
inequality~\eqref{op ineq}. Now $(J^{(n)})$ and $(J^{(n)*})$ are
sequences of contractions which converge strongly to $I$ and, for $v
\in \Dom Z$, $J^{(n)} v \Tends v$ in the graph norm of $Z$. Thus
$K^{(n)} \Tends Z$ on $\D_0$, $L^{(n)} = L' J^{(n)} \Tends L$ on
$\D_0$ (since $L'$ is $Z$-bounded) and $M^{(n)} = J^{(n)*} \ol{M}
\Tends M$ on $\D \uot \iD$. In other words $F^{(n)} \Tends F$
pointwise on $\D_0 \op (\D \uot \iD)$.
\end{proof}

\section{Stochastic Hille-Yosida}

In this section we obtain the stochastic generator of a strongly
continuous left contraction cocycle --- when it has one, an existence
theorem for the QSDE~\eqref{op LHP} is established, and some examples
are discussed. We also briefly describe the situation when $\noise$ is
separable and has a given orthonormal basis.

\subsection{Stochastic generator of a cocycle}

We first show that strongly continuous left contraction cocycles
satisfy a quantum stochastic differential equation under a minimal
condition for the equation to make sense --- namely that there is an
available dense domain for a coefficient operator to act on. It
amounts to a weak-differentiability condition (cf.\ \cite{AJL, FF8}).
Recall the notation~\eqref{FVT} and the identities~\eqref{DVT=DVD}
and~\eqref{FV,H,G}.

\begin{thm} \label{newswitch}
Let $V$ be a strongly continuous left contraction cocycle on $\init$
with noise dimension space $\noise$, let $\iT^\dagger$ and $\iT$ be
total subsets of $\noise$ containing $0$, and let $\iD = \Lin \iT$,
$\iD^\dagger = \Lin \iT^\dagger$ and $Z= H^V_{0,0}$. If $\D^{V,\iT}$
is dense in $\init$ then the following hold.
\begin{alist}
\item
For $F =F^{V,\iT}$, the process $V$ satisfies the operator
\tu{QSDE}~\eqref{op LHP} strongly on $\D^{V,0} \ot \e{0} +\D^{V,\iT}
\uot \Exp_\iD$.
\item
If $\D$ is a core for $Z$ contained in $\D^{V,\iT}$, then $V$ is the
unique contractive $\iD^\dagger$-weak solution of~\eqref{op LHP} on
$\D \uot \Exp_\iD$, for $F=F^{V,\iT}|_{\D \uot \iDhat}$.
\item
If $\D^{\Vtilde, \iT^\dagger}$ is also dense in $\init$ then $(F^{V,
\iT})^* \supset F^{\Vtilde, \iT^\dagger}$.
\end{alist}
\end{thm}

\begin{proof}
(a) Since $V$ is strongly measurable and contractive it suffices to
show that $V$ is a $\noise$-weak solution by the second remark after
Proposition~\ref{nec ineq}. But this follows from the semigroup
representation as follows. Let $u \in \init, f \in \step$ and $(v,g)
\in (\D^{V,0} \times \{0\}) \cup (\D^{V,\iD} \times \step_\iD)$. Then
by adaptedness and the semigroup representation~\eqref{semi decomp},
\[
\ip{u \e{f}}{V_t v\e{g}} = \ip{u}{P^{f(t_0),g(t_0)}_{t_1-t_0} \cdots
P^{f(t_n),g(t_n)}_{t-t_n} v} e^{\int^\infty_t \ip{f(s)}{g(s)} \, ds},
\]
and since $\D^{V,\iD} \subset \Dom \HVcd$ for all $c \in \noise$ and
$d \in \iD$ by Corollary~\ref{C1.4}, the (a.e.)\ derivative of this
with respect to $t$ is
\[
\bip{u}{P^{f(t_0),g(t_0)}_{t_1-t_0} \cdots P^{f(t_n),g(t_n)}_{t-t_n}
\bigl( H_{f(t_n),g(t_n)} - \ip{f(t)}{g(t)} \bigr) v} e^{\int^\infty_t
\ip{f(s)}{g(s)} \, ds},
\]
in other words $\ip{u\e{f}}{V_t E^{\wh{f}(t)} F E_{\wh{g}(t)} v\e{g}}$
by~\eqref{FV,H,G}. Thus $V$ satisfies~\eqref{op LHP} $\noise$-weakly
on $\D^{V,0} \ot \e{0} +\D^{V,\iD} \uot \Exp_\iD$.

(b) This follows from Theorem~\ref{uniqueness}.

(c) This follows from~\eqref{FV*}: $(F^{V,\iT})^* \supset (F^V)^*
\supset F^\Vtilde \supset F^{\Vtilde,\iT^\dagger}$.
\end{proof}

\begin{rem}
By Corollary~\ref{Mreg cocs} if $V$ is \emph{Markov}-regular then $F^V
\in B(\init \ot \khat)$, so $\D^{V,\noise} = \init$ and hence $V$
satisfies the operator QSDE on $\init \ot \Exp_\noise$ for this
bounded operator --- this is Theorem~6.7 of~\cite{father}.

The theorem also extends the main result of~\cite{AJL} to infinite
dimensional noise. Note that an application of the Banach-Steinhaus
Theorem is needed there in order to show that the form QSDE
coefficient is actually the component map of an operator. In infinite
dimensions the same argument again leads to a form QSDE for $V$,
however the Banach-Steinhaus Theorem does not help in this case. The
above result therefore also fills a gap in the proof of Theorem~2.4
of~\cite{FF8}.
\end{rem}

From this result and part (c) of Theorem~\ref{uniqueness} we may now
give necessary and sufficient conditions for a contraction cocycle to
satisfy a QSDE.

\begin{propn}
Let $V$ be a strongly continuous left contraction cocycle on $\init$
with noise dimension space $\noise$. Then the following are
equivalent.
\begin{rlist}
\item
$\DvT$ is dense in $\init$ for some total subset $\iT$ of $\noise$
containing $0$.
\item
$V$ strongly satisfies a \tu{QSDE} of the form~\eqref{op LHP} on some
domain of the form $\D_0 \ot \e{0} + \D \uot \Exp_\iD$.
\item
$V$ is a $\iT^\dagger$-$\iT$-solution of a form \tu{QSDE}~\eqref{form
LHP} on some domain of the form $\D \uot \Exp_\iT$.
\end{rlist}
\end{propn}

\begin{rem}
Thus if $V$ is a left contraction cocycle on $\init$ which satisfies a
QSDE of the type~\eqref{op LHP} on $\D^{V,0} \ot \e{0} + \D^{V,\iT}
\uot \Exp_\iD$, then
\[
F=F^{V,\iD}
\]
where $\iD = \Lin \iT$.
\end{rem}

\subsection{Coordinates}

Suppose that $\noise$ is separable with orthonormal basis $\eta =
(d_i)_{i \geqslant 1}$, and set $d_0 := 0$. Let $V$ be a strongly
continuous left contraction cocycle on $\init$ and suppose that $\D =
\bigcap_{\al,\be} \Dom \Gh{\al}{\be}$ is dense in $\init$, where
$\Gh{\al}{\be}$ denotes the generator $G^V_{c,d}$ for $c = d_\al$ and
$d = d_\be$. Then Theorem~\ref{newswitch} ensures that $V$ strongly
satisfies a Hudson-Parthasarathy equation
\[
dV_t = V_t \Fab \, d\Lambda^\be_\al (t), \quad V_0 =I,
\]
in which $[\Fab]_{\al, \be \geqslant 0}$ is the matrix of components
of an \emph{operator} $F \in \Op (\D \ot \iDhat)$ where $\iD = \Lin
\eta$ --- in other words the matrix is \dfn{semiregular} in the sense
that $\sum_{\al \geqslant 0} \norm{\Fab v}^2 < \infty$ for all $\be
\geqslant 0$ and $v \in \D$. Moreover the components are recovered
from the associated semigroup generators by the affine transformation
\[
\begin{aligned}
F^0_0 &= \Gh{0}{0} \\
F^i_0 &= \Gh{i}{0} -\Gh{0}{0} +\tfrac{1}{2}, \quad i \geqslant 1 \\
F^0_j &= \Gh{0}{j} -\Gh{0}{0} +\tfrac{1}{2}, \quad j \geqslant 1 \\
F^i_j &= \Gh{i}{j} -\Gh{i}{0} -\Gh{0}{j} +\Gh{0}{0} -\delta^i_j, \quad
i,j \geqslant 1,
\end{aligned}
\]
$\delta^i_j$ being the Kronecker delta.

\subsection{Cocycles from stochastic generators}

Our treatment of the existence question for~\eqref{op LHP} is founded 
on the following infinitesimal version of Theorem~\ref{semi approx}.

\begin{propn} \label{cocycle only}
Let $\Q_\iT = \{\Qcd: c,d \in \iT\}$ be a family of $C_0$-contraction
semigroups on $\init$, indexed by a total subset $\iT$ of $\noise$
containing $0$ and let $\Gcd$ denote the generator of the semigroup
$\Qcd$. Suppose that there is a sequence of strongly continuous left
contraction cocycles $(V^{(n)})$ on $\init$ and, for each $c,d \in
\iT$, a core $\D_{c,d}$ for $\Gcd$ such that
\begin{alist}
\item \label{Fn}
$\D_{c,d} \subset \D^{V^{(n)},\iT}$ for each $n \in \Nat$, and
\item \label{Ec}
$E^{\wh{c}} F^{(n)} E_{\wh{d}} - \chi(c,d) \Tends \Gcd$ pointwise on
$\D_{c,d}$, for all $c,d \in \iT$, where $F^{(n)} := F^{V^{(n)},\iT}$.
\end{alist}
Then there is a unique strongly continuous left contraction cocycle
$V$ whose associated semigroups include $\Q_\iT$. Moreover $V^{(n)}
\Tends V$ in the weak operator topology.
\end{propn}

\begin{proof}
We use the notation $\Qncd$ and $\Gcd^{(n)}$ for semigroups and
generators associated with the cocycle $V^{(n)}$.
Condition~\eqref{Fn} and Corollary~\ref{C1.4} imply that $\D_{c,d}
\uot \iDhat \subset \Dom F^{(n)}$, where $\iD = \Lin \iT$, so each
$F^{(n)}$ is densely defined, and also $E^{\wh{c}}F^{(n)} E_{\wh{d}} -
\chi(c,d) \subset G^{(n)}_{c,d}$ by~\eqref{FV,H,G}. Hence, by the
Trotter-Kato Theorem (\cite{Davies}, Corollary~3.18),
assumption~\eqref{Ec} implies that
\[
\lim_{n \Tends \infty} \sup_{t \in [0,T]} \norm{(Q^{(n)c,d}_t -
\Qcd_t) u} \Tends 0
\]
for all $c,d \in \iT$, $u \in \init$ and $T >0$. The result therefore
follows by Theorem~\ref{semi approx} and Lemma~\ref{cts coc}.
\end{proof}

\begin{rem}
This result is a stochastic generalisation of the Trotter-Kato
Theorem. In the usual version pointwise convergence of the generators
implies convergence of the sequence of semigroups in the strong
operator topology. However a similar strengthening of the conclusion
for cocycles is not possible --- as can be demonstrated using the
conditions for isometricity of cocycles given in terms of
conservativity of an associated quantum dynamical semigroup.
See~\cite{LWbedlewo} for details.
\end{rem}

\begin{thm} \label{gen approx}
Let $F \in \Op (\D \uot \iDhat)$ where $\D$ is a dense subspace of
$\init$ and $\iD = \Lin \iT$ for a total subset $\iT$ of $\noise$
containing $0$. Assume that
\begin{alist}
\item \label{bdd cont cond}
for each $c,d\in\iT$, $E^{\wh{c}}FE_{\wh{d}} - \chi (c,d)$ is a
pregenerator of a $C_0$-contraction semigroup $\Qcd$, and
\item \label{gens cge}
there is a sequence $(F^{(n)})$ in $B(\init \ot \khat)$ satisfying the
operator inequality~\eqref{op ineq}, such that, for all $c,d \in \iT$,
\[
E^{\wh{c}} F^{(n)} E_{\wh{d}} \Tends E^{\wh{c}} F E_{\wh{d}}
\ \text{ pointwise on } \D.
\]
\end{alist}
Then $F \subset F^{V,\iT}$ for a unique strongly continuous left
contraction cocycle $V$ on $\init$. Moreover, for all $c \in \noise$
and $d \in \iT$,
\begin{align}
G^V_{c,d} &= \ol{E^{\wh{c}} FE_{\wh{d}} - \chi (c,d)}, \text{ and }
\label{GV bar} \\
G^\Vtilde_{d,c} &\supset E^{\wh{d}} F^* E_{\wh{c}} - \chi (d,c).
\label{GV tilde}
\end{align}
\end{thm}

\begin{proof}
By Theorems~\ref{uniqueness} and~\ref{newswitch},
assumption~\eqref{bdd cont cond} (with $c=d=0$) implies uniqueness.
Let $\Gcd$ be the generator of $\Qcd$ and let $V^{(n)}$ be the
strongly continuous left contraction cocycle generated by $F^{(n)}$
(see the remark following Proposition~\ref{nec ineq}). Then the
hypotheses of Proposition~\ref{cocycle only} are satisfied with
$\D_{c,d} = \D$ for each $c,d \in \iT$. Let $V$ be the resulting
cocycle. Then $\Gcd^V = \Gcd$ so $E^{\wh{c}} F E_{\wh{d}} \subset
\Gcd^V + \chi (c,d)$ and therefore $\DvT \supset \D$. This gives
$F^{V,\iT} \supset F$ and so $G^V_{c,d} \supset \ol{E^{\wh{c}} F
E_{\wh{d}} - \chi (c,d)}$ for all $c \in \noise$ and $d \in \iD$,
by~\eqref{FV,H,G}. Now $\D$ is a core for $\Gcd^V$ when $c,d \in \iT$
so, by part~(b) of Proposition~\ref{generators compared}, it is also a
core when $c \in \noise$. The above inclusion is therefore an
equality. It remains only to verify the inclusion~\eqref{GV tilde},
but since $G^\Vtilde_{d,c} = (\Gcd^V)^*$ this follows by taking
adjoints.
\end{proof}

\begin{rems}
Under the conditions of the theorem, if also $\Dom F^* \supset
\D^\dagger \uot \wh{\iD^\dagger}$, for dense subspaces $\D^\dagger$
and $\iD^\dagger$ of $\init$ and $\noise$ respectively, then
$\D^{\Vtilde, \iD^\dagger}$ contains the dense subspace $\D^\dagger$
so $\Vtilde$ strongly satisfies the QSDE~\eqref{op LHP} with
coefficient $F^{\Vtilde, D^\dagger}$ on $\D^{\Vtilde, 0} \ot \e{0} +
\D^{\Vtilde, \iD^\dagger} \uot \Exp_{\iD^\dagger}$, and $F^* \supset
F^{\Vtilde, \iD^\dagger}$. In particular, $F^*$ satisfies the form
inequality~\eqref{form ineq} on $\D^{\Vtilde, 0} \op (\D^{\Vtilde,
\iD^\dagger} \uot \iD^\dagger)$.
\end{rems}

Our next result extends Fagnola's existence theorem (\cite{FF8}).
Whereas his proof requires separability of both of the Hilbert spaces
$\init$ and $\noise$, ours requires a strengthening of his condition
which amounts to $K$ being a pregenerator of a $C_0$-semigroup. This
difference in hypotheses reflects our difference of approach. Whereas
he approximates the solution \emph{process} by adapting Frigerio's
diagonalisation argument with the Arzel\`{a}-Ascoli Theorem to
cocycles constructed from bounded stochastic generators, we
approximate a sufficient number of the associated \emph{semigroup}
generators by exploiting the Trotter-Kato Theorem and this demands
stronger core requirements.

\begin{thm} \label{HYT}
Let $F \in \Op (\D \uot \iDhat)$, with block matrix form
$\begin{sbmatrix} K & M \\ L & C-I \end{sbmatrix}$, where $\D$ is a
dense subspace of $\init$ and $\iD = \Lin \iT$ for a total subset
$\iT$ of $\noise$ containing $0$. Suppose that
\begin{alist}
\item \label{nec form i}
$2 \re \ip{\xi}{F\xi} +\norm{\Delta F\xi}^2 \leqslant 0$ for all $\xi
\in \D \uot \iDhat$, and
\item \label{suff bi}
$K+ME_d -\frac{1}{2}\norm{d}^2$ is a pregenerator of a
$C_0$-semigroup, for each $d\in\iT$.
\end{alist}
Then $F\subset F^{V,\iT}$ for a unique strongly continuous left
contraction cocycle $V$ on $\init$, moreover~\eqref{GV bar}
and~\eqref{GV tilde} hold.
\end{thm}

\begin{proof}
By Proposition~\ref{ineq conseqs}, assumption~\eqref{nec form i}
implies that $\Gcd^0 := E^{\wh{c}} FE_{\wh{d}} -\chi (c,d)$ is
dissipative and $\Gcd := \ol{\Gcd^0}$ generates a $C_0$-contraction
semigroup if and only if $G_{b,d}$ does, for $b,c \in \noise$ and $d
\in \iD$. Since $K +ME_d - \frac{1}{2} \norm{d}^2 = G_{0,d}^0$, this
operator is dissipative for each $d \in \iD$ and
assumption~\eqref{suff bi} is equivalent to $G_{0,d}$ being a
$C_0$-semigroup generator for each $d \in \iT$. Therefore $\Gcd$ is
such a generator for each $c \in \noise$ and $d \in \iT$. In view of
Proposition~\ref{more ineq conseqs}, Theorem~\ref{gen approx}
therefore applies.
\end{proof}

\begin{cor}
If condition~\eqref{gens cge} of \tu{Theorem~\ref{HYT}} is replaced by
\begin{rlist}
\item[(bi)]
$\Ran (\lambda I -K)$ is dense in $\init$ for some $\lambda >0$, and
\item[(bii)]
$ME_d$ is $K$-bounded, for each $d \in \iT$,
\end{rlist}
then the conclusion of the theorem holds, moreover
\[
\Dom G^V_{e,d} \supset \Dom \ol{K} \text{ for all } e \in \noise
\text{ and } d \in \iD,
\]
with equality when $d$ is a sufficiently small multiple of an element
of $\iT$.
\end{cor}

\begin{proof}
Since $K$ is dissipative~(bi) is equivalent to $\ol{K}$ being a
$C_0$-contraction semigroup generator, by the Lumer-Phillips Theorem
(\cite{Davies}, Theorem~2.25). For $\al > 0$ let $\iT^\alpha = 
\bigl\{ (\al + \lambda_d)^{-1} d: d \in \iT \bigr\}$, where
$\lambda_d$ is the relative bound of $ME_d$ with respect to $K$, thus
$\Lin \iT^\alpha =\iD$ and $ME_d$ has relative bound less than one for
each $d \in \iT^\alpha$. Using Proposition~\ref{ineq conseqs} once
more this means that, for each $e \in \noise$ and $d \in \iT^\alpha$,
$E^{\wh{e}} F E_{\wh{d}} -\chi (e,d)$ is a relatively bounded
perturbation of $K$ with relative bound less than one. Therefore, by
Gustafson's Theorem, its closure has the same domain as $\ol{K}$ and
is a $C_0$-contraction semigroup generator, so~(b) of
Theorem~\ref{HYT} holds (with $\iT^\alpha$ in place of $\iT$) and the
theorem applies. We have
\[
\Dom G^V_{e,d} = \Dom \ol{K} \text{ for } e \in \noise, d \in
\iT^\alpha \text{ and } \alpha > 0.
\]
The proof therefore follows by Corollary~\ref{C1.4} since $\iD = \Lin
T^\al$.
\end{proof}

\begin{rem}
By choosing suitable functions $\mu$ and $\lambda$ in the example
below, it is possible to find an operator $F \in \Op (\D \uot \iDhat)$
satisfying the conditions of Theorem~\ref{HYT} but whose coefficients
$M_d$ are not $K$-bounded, so that $F$ is not covered by the corollary
above.
\end{rem}

Dualising we obtain alternative conditions.

\begin{cor}
Let $F \in \Op (\D \uot \iDhat)$ and $F^\dagger \in \Op (\D^\dagger
\uot \wh{\iD^\dagger})$ be densely defined operators on $\init \ot
\khat$ with block matrix forms $\begin{sbmatrix} K & M \\ L & C-I
\end{sbmatrix}$ and $\Bigl[ \begin{smallmatrix} K^\dagger & L^\dagger
\\ M^\dagger & C^\dagger -I \end{smallmatrix} \Bigr]$ respectively,
satisfying $F^* \supset F^\dagger$, where $\iD = \Lin \iT$ for a total
subset $\iT$ of $\noise$ containing $0$. Then the conclusions of
\tu{Theorem~\ref{HYT}} hold under the conditions
\begin{alist}
\item
$F$ and $F^\dagger$ satisfy the form inequality~\eqref{form ineq},
\item[(bi)]
$K^\dagger$ is a pregenerator of a $C_0$-semigroup on $\init$, and
\item[(bii)]
$\D$ is a core for the operator $(K^\dagger +E^d M^\dagger
-\frac{1}{2} \norm{d})^*$, for each $d \in \iT$.
\end{alist}
\end{cor}

\begin{proof}
In view of assumption~(bi), Proposition~\ref{ineq conseqs} applied to
$F^\dagger$ shows that $(K^\dagger +E^dM^\dagger -\frac{1}{2}
\norm{d}^2)$ is a pregenerator of a $C_0$-contraction
semigroup. Assumption (bii) therefore implies that its closure is
$(K +ME_d -\frac{1}{2} \norm{d}^2)^*$, thus $(K +ME_d -\frac{1}{2}
\norm{d}^2)$ is a pregenerator of a $C_0$-contraction semigroup and so
Theorem~\ref{HYT} applies.
\end{proof}

Many examples are covered by the following consequence of
Theorem~\ref{HYT}, with $C$ typically being unitary.

\begin{thm}
Let $H$ be a closed symmetric operator on $\init$, $L$ a closed
operator $\init \to \init \ot \noise$, $C$ a contraction operator on
$\init \ot \noise$ and $\iT$ a total subset of $\noise$ containing
$0$, such that $\D := \Dom H\cap \Dom L^*L \cap \bigcap_{d \in \iT}
\Dom L^* CE_d$ is dense in $\init$, and let $F = \begin{sbmatrix} K &
-L^*C \\ L & C-I \end{sbmatrix} \big|_{\D \uot \iDhat}$ where $K = iH
-\tfrac{1}{2} L^*L$ and $\iD = \Lin \iT$. Then the following hold.
\begin{alist}
\item
\begin{alist}
\item
$F$ satisfies~\eqref{form ineq}, with equality if and only if $C$ is
isometric.
\item
If there are constants $\gamma_d >0$ \tu{(}$d \in \iT$\tu{)} such that
\[
(\gamma_d I+\tfrac{1}{2} L^* L+L^*CE_d -iH)\D \text{ is dense in }
\init.
\]
then $F$ generates a strongly continuous left contraction cocycle on
$\init$.
\end{alist}
\item
Suppose that $\D^\dagger := \Dom H^* \cap \Dom L^* L \cap \bigcap_{d
\in \iT^\dagger} \Dom L^* E_d$ is dense in $\init$, and let
$\iD^\dagger = \Lin \iT^\dagger$ for another total subset of $\noise$
containing $0$.
\begin{alist}
\item
$F^*$ satisfies~\eqref{form ineq} on $\D^\dagger \uot
\wh{\iD^\dagger}$, with equality if and only if $C$ is coisometric on
$\init \ot \noise$. 
\item
If there are constants $\gamma_d >0$ \tu{(}$d \in \iT^\dagger$\tu{)}
such that
\[
(\gamma_d I +\tfrac{1}{2} L^* L - L^* E_d +iH^*) \D^\dagger \text{ is
dense in } \init
\]
then $F^\dagger := F^* \big|_{\D^\dagger \uot \wh{\iD^\dagger}}$
generates a strongly continuous left contraction cocycle.
\end{alist}
\end{alist}
\end{thm}

Examples in which $C=I$ and $H=0$ have arisen recently in the problem
of constructing stochastic dilations of quantum Markov
semigroups~\cite{sym CP}. In this case it suffices for $\D$ to be a
core for the positive selfadjoint operator $L^*L$, and for $L^* E_d \ 
(d \in \iT)$ to be relatively bounded with respect to $L^*L$.

Let $\D$ be the linear span of the standard orthonormal basis of
$\init := l^2 (\Int_+)$, let
\[
F = \begin{bmatrix} \nu(N) & W^* \ol{\lambda} (N) \\ -\lambda(N) W & 0
\end{bmatrix}, \text{ where } \nu(n) = i\mu(n) -\tfrac{1}{2}
|\lambda|^2(n+1),
\]
where $W$ and $N$ denote respectively the isometric right shift on
$\init$ and the number operator on $\init$, and $\lambda: \Int_+ \To
\Comp$ and $\mu: \Int_+ \To \Real$ are arbitrary functions. Then $F^*
+F+F^* \Delta F$ and $F+F^* +F\Delta F^*$ both vanish on $\D \op \D$
and Theorem~\ref{HYT} applies. Models of this type arise in the study
of inverse harmonic oscillators interacting with a heat bath in the
singular coupling limit (\cite{Wal}). Conditions on the pair
$(\lambda, \mu)$ which ensure isometry/unitarity of the resulting
contraction cocycle are investigated in~\cite{RHP}, from the point of
view of the right equation $dU_t =(F^* \ot I_\Fock) \wh{U}_t
d\Lambda_t$.
 
Classical birth and death processes have been constructed using
quantum stochastic calculus (\cite{FFBandD}, \cite{RHP}). These are
similarly covered by the above theorem, this time working with the
Hilbert space $l^2(\Int)$ and two dimensional quantum noise. These
examples and others are treated in detail in~\cite{LWbedlewo}.

\noindent
\emph{ACKNOWLEDGEMENTS}.
We are grateful to Nick Weatherall for drawing our attention to
Dixmier's counterexamples. SJW acknowledeges financial support for
this work from two EU TMR Networks (HPRN-CT-2002-00279 and
HPRN-CT-2002-00280).

\enlargethispage{17mm}

\end{document}